\documentclass{article}
\usepackage[utf8]{inputenc}
\usepackage{graphicx}

\title{Towards an untyped proof of Con(NF)}
\author{Zuhair Al-Johar}
\date{July 2021}

\begin{document}

\maketitle

\section{Introduction}

The idea of this approach towards proving the consistency of Quine's New Foundations set theory is to go in a completely untyped manner. So no contemplation about types is utilized here. All conceptualization pivots around proving a handful many axioms that are completely untyped and framed in the usual language of set theory, and proven to be equivalent to NF.  Here, it'll be shown that if we assume the consistency of ZF plus an automorphism and an external bijection with suitable internalization of subsets of its domain and codomain, then NF would be interpreted in this system using a modification of Boffa construction models.

\section {The ambient theory}

To the language of ZF we add two primitive unary functions  $j,f$ the first is an automorphism over the universe, the second is a partial function. The axioms are: \bigskip

\noindent
\textbf{ Axioms of $\sf ZF$ } but Separation and collection do not use symbols $f,j$ unless in parameters.\bigskip

\noindent

\noindent
\textbf{Injectivity:} $\forall x \, \forall y: f(x)=f(y) \to x=y$

\noindent
\textbf {Automorphism:} $\forall x \, \forall y: x \in y \iff j(x) \in j(y)$

\noindent
\textbf {Amenability:} $\forall a : \exists b (b=\{j(x):x \in a\}) \land \exists d (d=\{x: j(x) \in a\})$

\noindent
\textbf {Movement} $\exists \alpha: \lim(\alpha) \land V_{j(\alpha)+1} \subset   V_\alpha \land \\ f `` V_\alpha = V_{j(\alpha)+1}  \\   \forall x \subseteq V_\alpha \, \exists y: y =f``x  \\   \forall x \subseteq \mathcal  (V_{j(\alpha)+1})^2 \ \exists y: y=\{\{z,u\}: \{f(z),f(u)\} \in x\} $ \bigskip

\noindent
Where: $f``k= \{f(x): x \in k\}$, and $x^2= \{\{a,b\}: a,b \in x\}$ \bigskip

In English: we have the rank $V_{\alpha+1}$ (for limit $\alpha$) moved by automorphism $j$ to a proper subset of $V_\alpha$, that is we have $V_{j(\alpha)+1} \subset V_\alpha$. At the same time we have the external function $f$ being bijective from $V_\alpha$ to $V_{j(\alpha)+1}$ and such that for every subset $x$ of $V_\alpha$ there is a set-image of $x$ under $f$ (that is $f``x$), and also for every subset $y$ of $V_{j(\alpha)+1}$ there is pre-image of $y$ under $f$ (that is $f^{-1}``y$) (see Lemma below), and moreover for every set of pairs of elements of $V_{j(\alpha)+1}$ there is a set of all pairs of pre-images (under $f$) of those paired elements.  

\section {Interpreting NF}
\label{Interpreting}

Using Boffa construction models we know that $(V_\alpha, \in^*)$ would be a model of NFU + Infinity, where $\in^*$ defined as: $$y \in^* x \iff y \in j^{-1}(x) \land x \in V_{j(\alpha)+1}$$

Now we seek to prove that all members of $V_\alpha$ can be re-coded in a such a manner that turns them all to be sets! And of course at the same time keeping all rules of stratification. We do this using the external function $f$, so we need to define a new membership $\in'$ as:$$y \in' x \iff y \in j^{-1}(f(x))  $$

The stratification rules are those axiomtized by the following five axioms.

\textbf{Complements} $\forall x \exists y: y= \{z: z \not \in x\}$

\textbf{Pairing}: $\forall a \forall b \exists y: y=\{a,b\}$

\textbf{Set Union}: $\forall x \exists y: y= \bigcup x$

\textbf{Unordered Composition}: $\forall R \forall S \exists X: X= R \circ^* S$

\textbf{Unordered intersection relation set}: $\exists X: X= \Pi^* $ \smallskip

Where:
$c \circ^* d = \{\{x,z\}: \, \exists y \, (\{x,y\} \in c \land \{y,z\} \in d )\}$ \smallskip

$ \Pi^* A = \{\{x,y\} \in A: x \cap y \neq \emptyset \}$ \bigskip

That this is a finite axiomatization of Stratified Comprehension is presented in my article[1]. We denote this system by $\sf Fin. SF$\bigskip

So the idea is to replace each symbol $\in$ in the above by $\in'$ and restrict all quantifiers to $V_\alpha$, and prove the replacing system consistent relative to the ambient theory. \bigskip

\textbf{Lemma}: $\forall x \subseteq V_{j(\alpha)+1} \, \exists y: y =\{z: f(z) \in x\}$ \bigskip:

Proof: for every $x \subseteq V_{j(\alpha)+1}$ there is a set $\mathcal P_1(x)$ of all singletons of its elements, so by ``Movement" well have the set $\{\{z\}: \{f(z)\} \in \mathcal P_1(x)\}  $, take the union of this set and this would be the set $y$ above. \bigskip

\textbf {Proposition}: $(V_\alpha, \in^*) \models \sf Fin.SF$

Proof: Since $(V_\alpha, \in^*) \models \sf NFU$, then $(V_\alpha, \in^*) \models \sf Fin.SF$, since $\sf Fin.SF \subset NFU$.\bigskip

\textbf{The proof of Extensionality:} Since every element of $V_\alpha$ would be sent by $f$ to an element of $V_{j(\alpha)+1}$, then all of them would code (through $j^{-1} f$) elements of $V_{\alpha+1}$, and since the latter is $\in$-extensional, then all elements of $V_\alpha$  would be $\in'$-extensional, this follows from the definition of $\in'$ and from $j^{-1}f$ being a bijection! \bigskip

\textbf{The proof of Complements:}  for every element $x \in V_\alpha$, we have $f(x) \in V_{j(\alpha)+1}$ that is the $j$ code of an element $k \in V_{\alpha+1}$, which of course has its complement  $k^c \in V_{\alpha+1}$, which has $j(k^c) \in V_{j(\alpha)+1}$, now take $f^{-1}(j(k^c)) \in V_\alpha$ and this would be the complementary set of $x$ under the re-defined membership relation $\in'$ \bigskip

\textbf{The proof of Pairing:}  for all $a,b \in V_\alpha$ the set $f^{-1}(j(\{a,b\})) \in V_\alpha$, and this is the $\in'$-pair of $a,b$ \bigskip

\textbf{The poof of Set union:}  Let $l \in V_\alpha$, take $f(l)$ this would be an element of $V_{j(\alpha)+1}$ and so it is the $j$ code of a subset $x$ of $V_\alpha$ now $x$ is the set of all $\in'$-elements of $l$ [definition of $\in'$]. Now take $f``x$, this would be a subset of $V_\alpha$, so it has a $j$ code $j(f``x)$, now we know that $
(V_\alpha, \in^*)$ satisfy set unions, so we have an $\in^*$ set union of $j(f``x)$, denote  that by $k=\bigcup^{\in^*} j(f``x)$, now take $f^{-1}(k)$ and that would be the $\in'$ set union of $l$\bigskip

\textbf{The proof of Unordered Composition:} The proof is generally similar to set unions. For any sets $x,y \in V_\alpha$ we take $f(x),f(y)$, now the sets $k=j^{-1}(f(x)), l=j^{-1}(f(x))$ would be the sets of all $\in'$-elements of $x,y$ respectively. Now take $f``k,f``l$ those would be subsets of $V_\alpha$ and so have $j$-codes $j(f``k), j(f``l)$ Now those sets would have an $\in^*$ unordered composition of them, call that $q$. Now  we reverse the process, that is we take $j^{-1}(q)$ then take its pre-image under $f$ that is we take $f^{-1}``j^{-1}(q)$ call this is $r$, then take $f^{-1} (j(r))$ and this would be the $\in'$ unordered composition of $x,y$ \bigskip

\textbf{The proof of Unordered Intersection Relation Set} We start from the $\in^*$-intersection relation set, denote that by $\Pi^{\in^*}$. That is known to exist. Take $k = j^{-1} \Pi^{\in^*}$, now we take the set $\{\{z,u\}: \{f(z),f(u)\} \in k\}$ this would be the set of all pairs of intersecting $\in'$ elements of $V_\alpha$, call it $I$, now let $X=f^{-1} ``j(I)$,  then take $f^{-1}(j(X))$ and this would be our $\in'$-unordered intersection relation set.

\section{Remark:}
A proof of consistency of the ambient theory is needed to complete the proof of Con(NF). The peculiar thing about this approach is that in addition to its un-typed nature, it shows that the interpreting function $f$ which establishes the full Extensionality of NF, need not be an internal function at all. What is used to be known before is that if we prove the existence of an \emph{internal} bijection between the Ur-elements and sets of NFU, then we get an interpretation of NF. However, this proof shows that this need not be the case, and that external bijections with suitable additional internalization features can do the job. 

\section {Boffa models without automorphisms}
Boffa had used a rank shifting automorphism $j$ over a model  of $\sf ZF$ and showed that this would interpret $\sf NFU$, it'll be shown here that $j$ need not be an automorphism. All what is needed is for $j$ to be a partial unary function with the following features. However, I'll use the symbol $f$ instead since it resmbles the one used above. \bigskip

For some limit ordinal $\alpha$, we have: \smallskip

\textbf{Restriction:} $f : V_{\alpha+1} \to V_\alpha $ \bigskip

\textbf{Injectivity:} $f(x)=f(y) \to x=y$ \bigskip

\textbf{Upward:} $\forall x \subseteq V_\alpha \, \exists y: y=\{\{z,u\}: \{f(z),f(u)\} \in x \}$ \bigskip

\textbf{Downward:} $\forall x \subseteq V_{\alpha+2}\, \exists y: y=\{\{f(z),f(u)\}:  \{z,u\} \in x\}$ \bigskip

Now we KNOW that the above system is consistent, simply take $f$ to be an external automorphism that moves $V_{\alpha+1}$ to $V_{f(\alpha)+1} \subset V_\alpha$ , and all of the above rules would follow. However, the rules above doesn't prove the $f$ is an automorphism (see below). So, the conditions depicted here are weaker than those of an automorphism! \bigskip

\textbf{NOTE:} although the proof here uses ranks $V_\alpha, V_{\alpha+1}, f `` V_{\alpha+1}, etc..$, yet there is no need for that. All of what's needed is for $f$ to be an external injection whose domain is the power set of its codomain, that fulfills downward and upward axioms. However, we'll continue this tradition to conform more with the traditional approach. \bigskip

Now we set to prove the finite axiomatization of $\sf SF$ given above. \bigskip

We take $\in^*$ to come from domain $V_\alpha$, and is defined in terms of $f$ as:
$$y \in^* x \iff y \in f^{-1}(x) \land x \in V_{f(\alpha)+1}$$

We seek to interpret NFU over $\langle V_\alpha, \in^* \rangle$: \bigskip

\noindent
The proofs of \emph{Complements, Boolean union and Singletons} are straightforward.\bigskip

\textbf{Lemma 1}: $\forall x \subseteq V_\alpha\, \exists y: y =\{z: f(z) \in x\} = f^{-1}``x$ \bigskip

Proof: for every $x \subseteq V_\alpha$ there is a set $\mathcal P_1(x)$ of all singletons of its elements, so by ``Upward" we'll have the set $\{\{z\}: \{f(z)\} \in \mathcal P_1(x)\}  $, take the union of this set and this would be the set $y$ above. \bigskip

\textbf{Lemma 2}: $\forall x \subseteq V_{\alpha+1} \, \exists y: y =\{f(z): z \in x\} = f``x$ \bigskip

Proof: for every $ x \subseteq V_{\alpha+1}$ there is a set $\mathcal P_1(x)$ of all singletons of its elements, so by ``Downward" we'll have the set $\{\{f(z)\}: \{z\} \in \mathcal P_1(x)\}  $, take the union of this set and this would be the set $y$ above. \bigskip

\textbf{The proof of set unions:} let $ x \in V_\alpha $, let $f^{-1}(x) \subseteq V_\alpha$, now by Lemma 1 we'll have the set $k=\{z: f(z) \in f^{-1}(x)\}$, so $\bigcup k$ would be a subset of $V_\alpha$ and so  $f (\bigcup k)$  is the needed $\in^*$-set union of $x$ \bigskip

\textbf{Proof of Unordered Composition:} for $x,y \in V_\alpha$, let $f^{-1}(x), f^{-1}(y) \in V_{\alpha+1}$, then we take the sets $f^{-1}``f^{-1}(x), f^{-1}``f^{-1}(y)$, take the unordered composition of them, let that be $K$, then we reverse the process so $f(f``K)$ is the $\in^*$-unordered composition of $x$ and $y$ \bigskip

\textbf{Proof of the Unordered Intersection Relation Set}: We start from the set $\Pi=\{\{a,b\}: a \subseteq V_\alpha, b \subseteq V_\alpha, a \cap b \neq \emptyset \}$, that is: the set of all pairs of intersecting subsets of $V_\alpha$. Now this is a subset of $V_{\alpha+2}$, so by downward we have the set $K=\{ \{f(a),f(b)\}: \{a,b\} \in \Pi \}$, then by Lemma 2 we have $f``K$, then $f(f``K)$ is the $\in^*$-unordered intersection relation set. \bigskip

\textbf{Proof that $f$ is not necessarily an automorphism}: take the transposition $g$ of $\emptyset$ and $1 (i.e.;\{\emptyset\})$, now let $f=j \circ g= j(g(x))$ where $j$ is an automorphism over the universe. Take $f \uparrow V_{\alpha+1}$, now this is an injective function from $V_{\alpha+1} $  to $V_\alpha$, and it's easy to prove that it fulfills upwards and downwards. For downward direction if $x \subseteq \{\{a,b\}: a,b \subseteq V_\alpha\}$, then we simply take the subset $k$ of $x$ whose elements are pairs that do not intersect with $\{0,1\}$, now we have the set $A=\{\{f(a),f(b)\}: \{a,b\} \in k\}=\{\{j(a),j(b)\}: \{a,b\} \in k\}  $, because $f(s)=j(s)$ if $s \not \in \{0,1\}$; now for the rest of $x$, i.e. the set $k^c \cap x$, take the set $ B=\{ \{g(a),g(b)\}: \{a,b\} \in j(k^c \cap x)\}$, this can be easily constructed even in Zermelo, take $A \cup B$ and this would be the $f$-downward set. The same argument can be applied for the opposite direction to prove $f$-upward set. So $f$ fulfills all of the above axioms, yet clearly $f$ is not an automorphism.

\section {Significance}
It is $j$ being an automorphism in Boffa's construction that enforced having more Ur-elements than Sets in the interpretation of NFU. This would not necessarily be the case for the above function, so in principle it might be possible that $f``V_{\alpha+1}$ (the set of all sets in the interpretation) be of the same size or even strictly larger in size than its complementary set with respect to $V_\alpha$ (the set of all Ur-elements in the interpretation), in which case Con(NF) would follow. So in theory, the door is still open for a proof of Con(NF) using this method. 

\section {An aside: A proof of NF in NFU}
Along this method, it can be shown that if to axioms of NFU we add the following axiom: $$\exists x: |\mathcal P(x)|=|x| $$; where $||$ is for cardinality defined after Frege. Then NF follows. \bigskip

Here $\mathcal P(x) = \{y \in Set: \forall z \in y (z \in x)\}$, and  $Set=\{y: y=\emptyset \lor \exists x \, (x \in y)\}$

Proof: any bijection $f$ that witness the equality of cardinality between $x$ and its power set, would fulfill all of the above four axioms, and since $\mathcal P(x)$ is fully extensional (no distinct elements of it have exactly the same members), then Extensionality would be fulfilled.
So let $f: x \to P(x)$, then define $$y \in^f x \iff y \in f(x)$$, take the domain of $\in^f$ to be $x$, and we get $(x,\in^f)$  modeling $\sf NF$  \bigskip

\textbf{ Ur is not necessarily of empty objects }

Take a surjection $f: V \to Set$ that sends every empty object other than the empty set to some fixed set $x$, and send otherwise all elements of $Set$ to elements of $Set$. Apply the argument in the above section and we'll have an interpretation of NFU in which all Ur-elements are co-extensional to $x$.

\section{References}

\noindent
[1] Al-Johar, Z.A., Short Axiomatization of Stratified Comprehension, pre-print 2020, arXiv:2009.03185v2 [math.LO]

\noindent
[2]Ehrenfeucht, A., and A. Mostowksi, “Models of axiomatic theories admitting automorphisms,” Fundamenta Mathematicae, vol. 43 (1956), pp. 50–68. Zbl 0073.00704. MR 0084456. 574

\noindent
[3] M. Randall Holmes. "The Usual Model Construction for NFU Preserves Information." Notre Dame J. Formal Logic 53 (4) 571 - 580, 2012. 

\noindent
https://doi.org/10.1215/00294527-1722764

\noindent
[4] Quine, W. v. O., ``New foundations for mathematical logic", American Mathematical Monthly, vol. 44 (1937), pp. 70-80.

\end{document}